\documentclass[12pt]{amsart}
\usepackage{cases}
\usepackage{txfonts}
\usepackage{latexsym}
\usepackage{amsthm}
\usepackage{mathrsfs}
\usepackage{amssymb, amsmath}

\def\opn#1#2{\def#1{\operatorname{#2}}} % to make operators
\opn\chara{char} \opn\length{\ell}
\opn\projdim{proj\,dim} \opn\injdim{inj\,dim} \opn\rank{rank}
\opn\depth{depth} \opn\grade{grade} \opn\height{height}
\opn\embdim{emb\,dim} \opn\codim{codim}

\opn\Tr{Tr} \opn\bigrank{big\,rank}
\opn\superheight{superheight}\opn\lcm{lcm}
\opn\trdeg{tr\,deg}%
\opn\reg{reg} \opn\lreg{lreg}
\opn\div{div} \opn\Div{Div} \opn\cl{cl} \opn\Cl{Cl}
\opn\Spec{Spec} \opn\Supp{Supp} \opn\supp{supp} \opn\Sing{Sing}
\opn\Ass{Ass}
\opn\Ann{Ann} \opn\Rad{Rad} \opn\Soc{Soc}
\opn\Ker{Ker} \opn\Coker{Coker} \opn\Im{Im} \opn\Hom{Hom}
\opn\Tor{Tor} \opn\Ext{Ext} \opn\End{End} \opn\Aut{Aut} \opn\id{id}

\opn\nat{nat}
\opn\pff{pf}%   \pf exists already
\opn\Pf{Pf} \opn\GL{GL} \opn\SL{SL} \opn\mod{mod} \opn\ord{ord}
\opn\aff{aff} \opn\con{conv} \opn\relint{relint} \opn\st{st}
\opn\lk{lk} \opn\cn{cn} \opn\core{core} \opn\vol{vol}
\opn\gr{gr}

\def\pot#1#2{#1[\kern-0.28ex[#2]\kern-0.28ex]}

\opn\dirlim{\underrightarrow{\lim}}
\opn\invlim{\underleftarrow{\lim}}
\newtheorem{theorem}{Theorem}

\newtheorem{lemma}[theorem]{Lemma}
\newtheorem{conjecture}[theorem]{Conjecture}
\def\qed{\ifhmode\textqed\fi
   \ifmmode\ifinner\quad\qedsymbol\else\dispqed\fi\fi}
\def\textqed{\unskip\nobreak\penalty50
    \hskip2em\hbox{}\nobreak\hfil\qedsymbol
    \parfillskip=0pt \finalhyphendemerits=0}
\def\dispqed{\rlap{\qquad\qedsymbol}}

\opn\ini{in} \opn\inm{inm} \opn\Sym{Sym}

\begin{document}
\title{The $k$-tuple jumping champions among consecutive primes}
\author{XiaoSheng Wu and ShaoJi Feng}
\date{}
\address {Academy of Mathematics and Systems Science, Chinese Academy of Sciences,
Beijing 100190, P. R.
China.}
\email {xswu@amss.ac.cn}
\address {Academy of Mathematics and Systems Science, Chinese Academy of Sciences,
Beijing 100190, P. R.
China.}
\email {fsj@amss.ac.cn}

\begin{abstract} For any real $x$ and any integer $k\ge1$, we say that a set $\mathcal{D}_{k}$ of $k$ distinct integers is a $k$-tuple jumping champion if it is the most common differences that occurs among $k+1$ consecutive primes less than or equal to $x$. For $k=1$, it's known as the jumping champion introduced by J. H. Conway. In 1999 A. Odlyzko, M. Rubinstein, and M. Wolf announced the Jumping Champion Conjecture that the jumping champions greater than 1 are 4 and the primorials 2, 6, 30, 210, 2310,.... They also made a weaker and possibly more accessible conjecture that any fixed prime $p$ divides all sufficiently large jumping champions. These two conjectures were proved by Goldston and Ledoan under the assumption of appropriate forms of the Hardy-Littlewood conjecture recently.  In the present paper we consider the situation for any $k\ge2$ and prove that any fixed prime $p$ divides every element of all sufficiently large $k$-tuple jumping champions under the assumption that the Hardy-Littlewood prime $k+1$-tuple conjecture holds uniformly for $\mathcal{D}_k\subset[2,\log^{k+1}x]$. With a stronger form of the Hardy-Littlewood conjecture, we also proved that, for any  sufficiently large $k$-tuple jumping champion, the $gcd$ of elements in it is square-free.

\end{abstract}
\maketitle

\section{Introduction}
The study of finding the most probable difference among consecutive primes has existed for a long time. The problem was proposed by H. Nelson \cite{Nelson} in the issue of the 1977-78 volume of the Journal of Recreational Mathematics, and he had supposed 6 is the most probable difference between consecutive primes. However, assuming the prime pair conjecture from G. H. Hardy and J. E. Littlewood \cite{Littlewood}, P. Erd\"{o}s and E. G. Straus \cite{Erd}, in 1980, showed that there is no most likely difference, since they found that the most likely difference grows as the considered number becomes larger.

It was due to J. H. Conway who invented the term jumping champion to refer to the most common gap between consecutive primes not exceeding $x$. For the $n$th prime $p_n$, the jumping champions are the values of integer $d$ for which the counting function
\begin{align}
   N(x,d)=\sum_{{p_n\le x}\atop{p_n-p_{n-1}=d}}1\notag
\end{align}
attains its maximum
\begin{align}
   N^*(x)=\max_d N(x,d).\notag
\end{align}

In 1999 Odlyzko, Rubinstein and Wolf \cite{Odlyzko} announced the following two hypothesis, which are known as the Jumping Champion Conjecture now.
\begin{conjecture}
\label{A1}
   The jumping champions greater than 1 are 4 and the primorials 2, 6, 30, 210, 2310,$\cdots$.
\end{conjecture}
\begin{conjecture}
\label{B1}
   The jumping champions tend to infinity. Furthermore, any fixed prime $p$ divides all sufficiently large jumping champions.
\end{conjecture}

It's obvious that Conjecture \ref{B1} is a weaker consequence of Conjecture \ref{A1}, and as already mentioned, the first assertion of Conjecture \ref{B1} was proved by Erd\"{o}s and Straus \cite{Erd}, under the assumption of the Hardy-Littlewood prime pair conjecture. Recently, Goldston and Ledoan \cite{Goldston} extended successfully Edor\"{o}s and Straus's method to give a complete proof of Conjecture \ref{B1} under the same assumption. Soon after, they also give a proof of Conjecture \ref{A1} by assuming a sufficiently strong form of the Hardy-Littlewood prime pair conjecture.

Motivated by the work of Goldston and Ledoan, we have been working on the problem what are the most probable differences among $k+1$ consecutive primes with any $k\ge1$.

Let $\mathcal{D}_k=\{d_1,d_2,\cdots,d_k\}$ be a set of $k$ distinct integers with $d_1<d_2<\cdots<d_k$. For the $n$th prime $p_n$, we define the $k$-tuple jumping champions are the sets of $\mathcal{D}_{k}$ for which the sum
\begin{align}
   N_k(x,\mathcal{D}_{k})=\sum_{{p_{n+k}\le x}\atop{p_{n+i}-p_n=d_i}}1\notag
\end{align}
attains its maximum
\begin{align}
   N_k^*(x)=\max_{\mathcal{D}_{k}}N_k(x,\mathcal{D}_{k}).\notag
\end{align}

In the present paper, we work on the $k$-tuple jumping champion, and our main result can be summarized as follows.
\begin{theorem}
\label{C1}
   Let $k$ be any given positive integer. Assume Conjecture \ref{E1}. The $gcd$ (greatest common divisor) of all elements in the $k$-tuple jumping champions tend to infinity. Furthermore, any fixed prime $p$ divides every element of all sufficiently large $k$-tuple jumping champions.
\end{theorem}
With a stronger form of the Hardy-Littlewood conjecture, we obtain a stronger result.
\begin{theorem}
\label{A2}
    Assume Conjecture \ref{B2}, the $gcd$ of any sufficiently large $k$-tuple jumping champion is square-free.
\end{theorem}

In the following, we will denote $\mathcal{D}_k=d*\mathcal{D}'_k$, where $d=(d_1,d_2,\cdots,d_k)$ is the $gcd$ of the elements in $\mathcal{D}_k$ and $\mathcal{D}'_k=\{d_1',d_2',\cdots,d_k'\}$ with $d_i=dd_i'$ for any $i\le k$. We announce here that $\epsilon$ always denotes an arbitrary small positive constant but may have different value according to the context.

\section{the hardy-littlewood prime $n$-tuple conjecture}
Let $\pi_n(x,\mathcal{D}_n)$ denote the number of positive integers $m\le x$ such that $m+d_1,m+d_2,\cdots m+d_n$ are all primes and $\nu_{\mathcal{D}_n}(p)$
represents the number of distinct residue classes modulo $p$ occupied by elements of $\mathcal{D}_n$. The $n$-tuple conjecture says
\begin{align}
   \pi_n(x,D_n)\sim\mathfrak{S}(\mathcal{D}_n)\int_2^x\frac{dt}{\log^{n}t}\notag
\end{align}
as $x\rightarrow\infty$, where
\begin{align}
   \mathfrak{S}(\mathcal{D}_n)=\prod_{p}\Big(1-\frac1p\Big)^{-n}\Big(1-\frac{\nu_{\mathcal{D}_n}(p)}{p}\Big),\notag
\end{align}
with $p$ runs through all the primes.

In the proof of Theorem \ref{C1}, we need the following conjecture.
\begin{conjecture}
\label{E1}
If $\mathfrak{S}(\{0\}\cup\mathcal{D}_{k})\neq0$, as $x\rightarrow\infty$
\begin{align}
   \pi_{k+1}(x,\{0\}\cup\mathcal{D}_{k})=\mathfrak{S}(\{0\}\cup\mathcal{D}_{k})\frac{x}{\log^{k+1}x}(1+o(1))\notag
\end{align}
uniformly for $\mathcal{D}_{k}\subset[2,\log^{k+1}x]$.
\end{conjecture}

It is reasonable to suppose that the Hardy-Littlewood conjecture will hold uniformly for any $\mathcal{D}_{k}\subset[2,x]$, but the range $[2,\log^{k+1}x]$ is enough for our proof.

To prove Theorem \ref{A2}, we need the following stronger form of the Hardy-Littlewood Conjecture.

\begin{conjecture}
\label{B2}
For $n=k+1, k+2$, if $\mathfrak{S}(\{0\}\cup\mathcal{D}_{n})\neq0$, as $x\rightarrow\infty$
\begin{align}
   \pi_{n}(x,\{0\}\cup\mathcal{D}_{n-1})=\mathfrak{S}(\{0\}\cup\mathcal{D}_{n-1})\frac{x}{\log^{n}x} \bigg(1+E_n\bigg)\notag
\end{align}
uniformly for $\mathcal{D}_{n}\subset[2,\log^{k+1}x]$, where
\begin{align}
   E_n=\left\{
\begin{array}
    {r@{\quad:\quad}l}
    o\bigg(\frac1{(\log\log x)^2}\bigg) & n=k+1; \\
    o(1)  & n=k+2 .
\end{array}
\right.\notag
\end{align}
\end{conjecture}
We also need the following well-known sieve bound, for $x$ sufficiently large,
\begin{align}
\label{a2}
   \pi_n(x,\mathcal{D}_{n})\le(2^nn!+\epsilon)\mathfrak{S}(\mathcal{D}_{n}) \frac{x}{\log^nx},
\end{align}
for $\mathfrak{S}(\mathcal{D}_{n})\neq0$, which was given by Halberstam and Richert's excellent monograph \cite{Halberstam}.
\section{Lemma}
To prove Theorem \ref{A2}, we need the following lemmas.
\begin{lemma}
\label{C2}
For any set $\mathcal {D}_k\subset[0,h]$, $H\le h$, we have
\begin{align}
   \sum_{{1\le d_0\le H}\atop{d_0\notin\mathcal{D}_k}}\mathfrak{S}(\mathcal {D}_k\cup\{d_0\})=\mathfrak{S}(\mathcal {D}_k)H\Big(1+O_k\big(\frac{h^{\epsilon}}{H^{1/2}}\big)\Big).\notag
\end{align}
\end{lemma}
This lemma is about the average of the singular series, and the study of this is interesting in itself. We will give the proof of this in the last section.
\begin{lemma}
\label{D2}
For any integer $k\ge1$, assume Conjectures \ref{B2}.  Let $\mathcal{D}_k$ be a set of $k$ distinct integers with $\mathfrak{S}(\{0\}\cup\mathcal{D}_{k})\neq0$.
\begin{description}
  \item[(i)] If $2\le d_{k}=o(\log x)$, then
\begin{align}
   N_k(x,\mathcal{D}_{k})=\mathfrak{S}(\{0\}\cup\mathcal{D}_{k})\frac{x}{\log^{k+1}x}\bigg\{1-\frac{d_{k}}{\log x}+o\bigg(\frac{d_{k}}{\log x}\bigg)+o\bigg(\frac1{(\log\log x)^2}\bigg)\bigg\}.\notag
\end{align}
  \item[(ii)] If $H\le d_{k}\le\log^{k+1} x$ for some $H$ with $\log x/\log\log x\le H= o(\log x)$, then
\begin{align}
   N_k(x,\mathcal{D}_{k})\le\mathfrak{S}(\{0\}\cup\mathcal{D}_{k})\frac{x}{\log^{k+1}x}\bigg\{1-\frac{H}{\log x}+o\bigg(\frac{H}{\log x}\bigg)\bigg\}.\notag
\end{align}
\end{description}
\end{lemma}
\proof By inclusion-exclusion we have, for any integer $I\ge0$ and any $1\le H\le d_{k}$, we have
\begin{align}
\label{b2}
   N_k(x,\mathcal{D}_{k})\ge\sum_{i=0}^{2I+1}(-1)^{i}\sum_{{0<m_1<\cdots<m_i<d_{k}}\atop{m_1,\cdots,m_i\notin \mathcal{D}_{k}}}\pi_{k+1+i}(x,\{0,m_1,\cdots,m_i\}\cup\mathcal{D}_{k}) ,
\end{align}
and
\begin{align}
\label{c2}
   N_k(x,\mathcal{D}_{k})\le\sum_{i=0}^{2I}(-1)^{i}\sum_{{0<m_1<\cdots<m_i<H}\atop{m_1,\cdots,m_i\notin \mathcal{D}_{k}}}\pi_{k+1+i}(x,\{0,m_1,\cdots,m_i\}\cup\mathcal{D}_{k}).
\end{align}
By Conjecture (\ref{B2}) and Lemma (\ref{C2}), we find, for sufficiently large $x$,
\begin{align}
\label{d2}
   \sum_{{0<m_1<H}\atop{m_1\notin \mathcal{D}_{k}}}\pi_{k+2}(x,\{0,m_1\}\cup\mathcal{D}_{k})
   =&\sum_{{0<m_1<H}\atop{m_1\notin \mathcal{D}_{k}}}\mathfrak{S}(\{0,m_1\}\cup\mathcal{D}_{k})\frac{x}{\log^{k+2}x} (1+o(1))\notag\\
   =&\mathfrak{S}(\{0\}\cup\mathcal{D}_{k})\frac{x}{\log^{k+1}x}\frac{H}{\log x}\bigg(1+O\bigg(\frac{d_{k}^\epsilon}{H^{1/2}}\bigg)+o(1)\bigg).
\end{align}
From (\ref{a2}) and Lemma \ref{C2}, we also have, for any $1\le H\le d_k$,
\begin{align}
\label{e2}
   \sum_{{0<m_1<m_2<H}\atop{m_1,m_2\notin \mathcal{D}_{k}}}\pi_{k+3}(x,\{0,m_1,m_2\}\cup\mathcal{D}_{k})&\ll \sum_{{0<m_1<m_2<H}\atop{m_1,m_2\notin \mathcal{D}_{k}}}\mathfrak{S}(\{0,m_1,m_2\}\cup\mathcal{D}_{k})\frac{x}{\log^{k+3}x}\notag\\
   &\ll\mathfrak{S}(\{0\}\cup\mathcal{D}_{k})\frac{x}{\log^{k+1}x}\bigg(\frac{H}{\log x}\bigg)^2\bigg(1+O\Big(\frac{d_k^\epsilon}{H^{1/2}}\Big)\bigg)^2.
\end{align}
In the process to obtain (\ref{d2}) and (\ref{e2}),  we ignore the terms with $\mathfrak{S}(\{0,m_1,\cdots,m_i\})=0$, since these terms have $\pi_{k+i}(x,\{0,m_1,\cdots,m_i\}\cup\mathcal{D})=0$ or 1 and contribute $\ll H^i$, which is absorbed in the error term.

Then employ (\ref{d2}) and (\ref{e2}) into (\ref{c2}) with $I=1$, we have
\begin{align}
  N_k(x,\mathcal{D}_{k})&\le\mathfrak{S}(\{0\}\cup\mathcal{D}_{k})\frac{x}{\log^{k+1}x}\bigg\{1-\frac{H}{\log x}\bigg(1+O\Big(\frac{d_{k}^\epsilon}{H^{1/2}}\Big)+o(1)\bigg)\notag\\
  &\ \ \ +O\bigg(\frac{H}{\log x}\bigg)^2\bigg(1+O\Big(\frac{d_k^\epsilon}{H^{1/2}}\Big)\bigg)^2\bigg\}\notag\\
  &\le\mathfrak{S}(\{0\}\cup\mathcal{D}_{k})\frac{x}{\log^{k+1}x}\bigg\{1-\frac{H}{\log x}+o\bigg(\frac{H}{\log x}\bigg)\bigg\}\notag
\end{align}
for any $H$ with $\log x/\log\log x\le H= o(\log x)$ and $H\le d_k\le\log^kx$ since $\epsilon$ can be chosen arbitrary small. Hence,we proved part (ii) of the lemma. To prove part (i), we set $H=d_k$ in (\ref{d2}) and (\ref{e2}). Since $2\le d_{k}=o(\log x)$, then part (i) follows by substituting (\ref{d2}) and (\ref{e2}) into (\ref{b2}) and (\ref{c2}) with $I=1$.

\section{proof of Theorem \ref{C1}}
We will only give the proof of the theorem for $k\ge2$, since the situation of $k=1$ has been proved by Goldston and Ledoan \cite{Goldston}.

It's not difficult for us to see
\begin{align}
\label{a}
   \pi_{k+1}(x,\{0\}\cup\mathcal{D}_{k})-\sum_{{d'<d_{k}}\atop{d'\notin\mathcal{D}_{k}}} \pi_{k+2}(x,\{0,d'\}\cup\mathcal{D}_{k})\le N_k(x,\mathcal{D}_{k})\le\pi_{k+1}(x,\{0\}\cup\mathcal{D}_{k}).
\end{align}
Therefore, by inequality (\ref{a2}), $\mathfrak{S}(\mathcal{D}_{n})\le d_{n}^\epsilon$ (the proof of this inequality is the same to section 4 of \cite{Goldston}) and $\pi_n(x,\mathcal{D}_{n})=0$ or 1 for $\mathfrak{S}(\mathcal{D}_{n})=0$, it follows that
\begin{align}
   \sum_{{d'<d_{k}}\atop{d'\notin\mathcal{D}_{k}}} \pi_{k+2}(x,\{0,d'\}\cup\mathcal{D}_{k})\ll d_{k}^{1+\epsilon}\frac{x}{\log^{k+2}x}.\notag
\end{align}
Hence, with the condition given by Theorem \ref{C1}, we have
\begin{align}
\label{b}
   N_k(x,\mathcal{D}_{k})=\mathfrak{S}(\{0\}\cup\mathcal{D}_{k})\frac{x}{\log^{k+1}x}(1+o(1)), \ \ \ \text{uniformly for} \ \ 2\le d\le(\log x)^{1-\epsilon}
\end{align}
and
\begin{align}
\label{c}
   N_k(x,\mathcal{D}_{k})\le\mathfrak{S}(\{0\}\cup\mathcal{D}_{k})\frac{x}{\log^{k+1}x}(1+o(1)), \ \ \ \text{uniformly for} \ \ 2\le d\le\log^{k+1}x.
\end{align}

In the following, we define
\begin{align}
   \mathcal{P}_n:=2\cdot3\cdot5\cdots p_n\notag
\end{align}
to denote the $n$-th term in the sequence of primorials and use $\lfloor y\rfloor$ to be the largest primorial not greater than $y$. Let ${\mathcal{K}}= \{1,2,\cdots, k\}$,  from (\ref{b}), it follows that
\begin{align}
\label{d}
   \mathfrak{S}(\{0\}\cup\lfloor(\log x)^{1/2}\rfloor*\mathcal{K})\frac{x}{\log^{k+1}x}(1-o(1))\le\max_{2\le d_{k}\le (\log x)^{1-\epsilon}}N_k(x,\mathcal{D}_k)\le N_k^*(x).
\end{align}
Here the choice of $\mathcal{K}$ is insignificant. In fact, it can be replaced by any bounded set of k coprime positive integers. On the other hand,
\begin{align}
   N_k(x,\mathcal{D}_{k})\le\sum_{{p_n\le x}\atop{p_n-p_{n-k}\ge d_{k}}}1\le\sum_{{p_n\le x}\atop{p_n-p_{n-k}\ge d_{k}}}\frac{p_n-p_{n-k}}{d_{k}}\le\frac{kx}{d_{k}},\notag
\end{align}
we have
\begin{align}
   N_k(x,\mathcal{D}_{k})\le\frac{kx}{\log^{k+1}x}, \ \ \ \texttt{for} \ \ d_k\ge\log^{k+1}x.\notag
\end{align}
However, from (\ref{d}) we have
\begin{align}
    N_k^*(x)\ge\mathfrak{S}(\{0\}\cup\lfloor(\log x)^{1/2}\rfloor*\mathcal{K})\frac{x}{\log^{k+1}x}(1-o(1)),\notag
\end{align}
while
\begin{align}
   \mathfrak{S}(\{0\}\cup\lfloor(\log x)^{1/2}\rfloor*\mathcal{K})&\ge\prod_{p\le(1/2-\epsilon)\log\log x}\Big(1-\frac1p\Big)^{-k}\prod_{p\ge(1/2-\epsilon)\log\log x}\Big(1-\frac1p\Big)^{-(k+1)}\Big(1-\frac {k+1}p\Big)\notag\\
   &\gg\prod_{p\le(1/2-\epsilon)\log\log x}\Big(1-\frac1p\Big)^{-k}\notag\\
   &\gg\exp\bigg(k\sum_{p\le(1/2-\epsilon)\log\log x}\frac1p+o(1)\bigg)\notag\\
   &\gg(\log\log\log x)^{k}\notag
\end{align}
by an application of Merten's formula (see Ingham's tract \cite{Ingham}, Theorem 7, Formula (23), p. 22)
\begin{align}
   \sum_{p\le x}\frac1p=\log\log x+B+O\Big(\frac 1{\log x}\Big), \ \ \ \texttt{as} \ \ x\rightarrow\infty,\notag
\end{align}
with $B$ is a constant. Hence, as $x$ sufficiently large, if $\mathcal{D}_{k}$ is a $k$-tuple jumping champion, then $d_{k}\le\log^{k+1}x$.

For $m\ge3$, let $\mathcal{D}_m=\{d_1,d_2,\cdots,d_m\}$ be a set of $m$ distinct integers with $d_1<d_2<\cdots<d_m$ and $d_{m}\le(\log x)^A$ for any given positive $A>1$. Since
\begin{align}
   \mathcal{D}_{m}=d*\mathcal{D}'_{m}\notag
\end{align}
with $d$ is the $gcd$ of all elements in $\mathcal{D}_m$, it's obvious that $d'_{m}<(\log x)^A$ and
\begin{align}
\label{a1}
  2\le\nu_{\mathcal{D}'_{m}}(p)=\nu_{\mathcal{D}_{m}}(p)\le k
\end{align}
for $p\nmid d$. Let
\begin{align}
   \Delta_{\mathcal{D}_m'}=\prod_{j< i}(d_i'-d_j')\notag
\end{align}
and $\omega(n)$ be the number of prime factors (not the number different prime factors) contained in positive integer $n$. Then from the well known fact, for sufficiently large integer $n$,
\begin{align}
   \omega(n)\le(1+\epsilon)\log n/\log\log n,\notag
\end{align}
we see that, for sufficiently large $x$
\begin{align}
\label{c1}
   \omega(\Delta_{\mathcal{D}_m'})\le\frac{Am(m-1)}{2}(1+\epsilon)\log\log x/\log\log\log x.
\end{align}
Furthermore, if $\nu_{\mathcal{D}'_{m}}(p)<m$, it means that $p\mid\Delta_{\mathcal{D}_m'}$. We see that the number of such $p$ with $\nu_{\mathcal{D}'_{m}}(p)<m$ is not more than $\frac{Am(m-1)}{2}(1+\epsilon)\log\log x/\log\log\log x$ for sufficiently large $x$. Then, from the definition of $\mathfrak{S}(\mathcal{D}_{m})$ and (\ref{a1}), we have
\begin{align}
\label{b1}
   \mathfrak{S}(\mathcal{D}_{m})&=\prod_{p}(1-\frac1p)^{-m}\prod_{p\mid d}\Big(1-\frac1p\Big)\prod_{{p\nmid d}\atop{p\mid\Delta_{\mathcal{D}_m'}}}\Big(1-\frac{\nu_{\mathcal{D}'_{m}}(p)}{p}\Big)\prod_{{p\nmid d}\atop{p\nmid\Delta_{\mathcal{D}_m'}}}\Big(1-\frac{m}{p}\Big)\notag\\
   &\le\prod_{p}(1-\frac1p)^{-m}\prod_{p\mid d}\Big(1-\frac1p\Big)\prod_{{p\nmid d}\atop{p\mid\Delta_{\mathcal{D}_m'}}}\Big(1-\frac{2}{p}\Big)\prod_{{p\nmid d}\atop{p\nmid\Delta_{\mathcal{D}_m'}}}\Big(1-\frac{m}{p}\Big).
\end{align}
Let $d'$ be the greatest square-free factor of $d$. It's obvious that $\omega(\lfloor d\rfloor)\ge\omega(d')$. Since the combination of the last three products in the last expression of (\ref{b1}) takes over all primes, we have
\begin{align}
\label{k}
   \mathfrak{S}(\mathcal{D}_{m})&\le\prod_{p}(1-\frac1p)^{-m}\prod_{p\mid d'}\Big(1-\frac1p\Big)\prod_{{p\nmid d'}\atop{p\mid\Delta_{\mathcal{D}_m'}}}\Big(1-\frac{2}{p}\Big)\prod_{{p\nmid d'}\atop{p\nmid\Delta_{\mathcal{D}_m'}}}\Big(1-\frac{m}{p}\Big)\notag\\
   &\le\prod_{p}(1-\frac1p)^{-m}\prod_{p\le p_{\omega(d')}}\Big(1-\frac1p\Big)\prod_{{p> p_{\omega(d')}}\atop{p\mid\Delta_{\mathcal{D}_m'}}}\Big(1-\frac{2}{p}\Big)\prod_{{p> p_{\omega(d')}}\atop{p\nmid\Delta_{\mathcal{D}_m'}}}\Big(1-\frac{m}{p}\Big)\notag\\
   &\le\prod_{p}(1-\frac1p)^{-m}\prod_{p\mid \lfloor d\rfloor}\Big(1-\frac1p\Big)\prod_{{p\nmid \lfloor d\rfloor}\atop{p\mid\Delta_{\mathcal{D}_m'}}}\Big(1-\frac{2}{p}\Big)\prod_{{p\nmid \lfloor d\rfloor}\atop{p\nmid\Delta_{\mathcal{D}_m'}}}\Big(1-\frac{m}{p}\Big).
\end{align}
Here, the sceond inequality in (\ref{k}) holds because that we may interchange every prime greater than $p_{\omega(d')}$ in the second product with a prime less than $p_{\omega(d')+1}$ in the last two products with an increase of the value to the formula. These interchanges can be made also because the fact that the combination of the last three products in the formula takes over all primes. The last inequality in (\ref{k}) is a result of the fact $\omega(d')\le\omega(\lfloor d\rfloor)$.
Let $\mathcal{M}=\{1,2,\cdots,m-1\}$, from (\ref{k}) and the inequality $\lfloor d\rfloor\le\lfloor\log^Ax\rfloor$, we may have
\begin{align}
\label{e1}
   \frac{\mathfrak{S}(\mathcal{D}_{m})}{\mathfrak{S}(\{0\}\cup\lfloor\log^Ax\rfloor*\mathcal{M})} &\le\frac{\prod_{p\mid \lfloor d\rfloor}\Big(1-\frac1p\Big)\prod_{{p\nmid \lfloor d\rfloor}\atop{p\mid\Delta_{\mathcal{D}_m'}}}\Big(1-\frac{2}{p}\Big)\prod_{{p\nmid \lfloor d\rfloor}\atop{p\nmid\Delta_{\mathcal{D}_m'}}}\Big(1-\frac{m}{p}\Big)}{\prod_{p\mid \lfloor\log^Ax\rfloor}\Big(1-\frac1p\Big)\prod_{p\nmid \lfloor\log^Ax\rfloor}\Big(1-\frac{m}{p}\Big)}\notag\\
   &\le\frac{\prod_{{p\nmid \lfloor \log^Ax\rfloor}\atop{p\mid\Delta_{\mathcal{D}_m'}}}\Big(1-\frac{2}{p}\Big)}{\prod_{{p\nmid \lfloor\log^Ax\rfloor}\atop{p\mid\Delta_{\mathcal{D}_m'}}}\Big(1-\frac{m}{p}\Big)}\le\prod_{{p\nmid \lfloor\log^Ax\rfloor}\atop{p\mid\Delta_{\mathcal{D}_m'}}}\frac{p-2}{p-m}.
\end{align}
 Then, by the prime number theorem and (\ref{c1}), we have
\begin{align}
\label{e}
   \frac{\mathfrak{S}(\mathcal{D}_{m})}{\mathfrak{S}(\{0\}\cup\lfloor\log^Ax\rfloor*\mathcal{M})}\le\prod_{{p\nmid \lfloor\log^Ax\rfloor}\atop{p\mid\Delta_{\mathcal{D}_m'}}}\frac{p-2}{p-m}
   &\le \prod_{\omega([\log^Ax])<p\leq \omega([\log^Ax])+\omega(\Delta_{D_m'})} \frac{p-2}{p-m}\notag\\
   &\le\prod_{(A-\epsilon)\log\log x\le p\le (A+Am(m-1)/2+\epsilon)\log\log x}\frac{p-2}{p-m}.
\end{align}
By the Meterns formula, we obtain the last expression in (\ref{e}) is
\begin{align}
   &\le\exp\bigg(\sum_{(A-1)\log\log x\le p\le (Am(m-1)/2+A+1)\log\log x}\log\Big(1+\frac{m-2}{p-m}\Big)\bigg)\notag\\
   &\le\exp\bigg(\sum_{(A-1)\log\log x\le p\le (Am(m-1)/2+A+1)\log\log x}\frac{m-2}{p}+O\bigg(\sum_{n\ge(A-1)\log\log x}\frac1{n^2}\bigg)\bigg)\notag\\
   &\le\exp\Big(\log^{m-2}\Big(\frac{\log\log\log x+\log(Am(m-1)/2+A+1)}{\log\log\log x+\log(A-1)}\Big)+O\Big(\frac1{\log\log\log x}\Big)\Big)\notag\\
   &\le1+O\Big(\frac1{\log\log\log x}\Big).\notag
\end{align}
Thus we have
\begin{align}
\label{f}
   \frac{\mathfrak{S}(\mathcal{D}_{m})}{\mathfrak{S}(\{0\}\cup\lfloor\log^Ax\rfloor*\mathcal{M})} \le1+O\Big(\frac1{\log\log\log x}\Big),
\end{align}
for any $A>1$ given.

From now on, we use $\mathcal{D}^*_{k}$ to denote a $k$-tuple jumping champion. Let $p^*<\log x$ is a given prime that $p^*\mid\lfloor\log^{k+1}x\rfloor$ but $p^*\nmid d^*$, it's obvious that $p^*d_{k}\le\log^{k+2}x$. Then using (\ref{f}) with
\begin{align}
   \mathcal{D}_{m}=\{0\}\cup p^**\mathcal{D}^*_{k}\notag
\end{align}
and
\begin{align}
   A=k+2\notag
\end{align}
 we can see
\begin{align}
\label{g}
   \mathfrak{S}(\{0\}\cup\mathcal{D}^*_{k})\Big(1+\frac{\nu_{p^*}(\{0\}\cup\mathcal{D}'^*_{k})-1} {p^*-\nu_{p^*}(\{0\}\cup\mathcal{D}'^*_{k})}\Big)&=\mathfrak{S}(\{0\}\cup p^**\mathcal{D}^*_{k})\notag\\
   &\le\mathfrak{S}(\{0\}\cup\lfloor\log^{k+2}x\rfloor*\mathcal{K}) \Big(1+O\Big(\frac1{\log\log\log x}\Big)\Big).
\end{align}
Here, $\nu_{p^*}(\{0\}\cup\mathcal{D}'^*_{k})<p^*$. This is because that $\pi_{k+1}(x, \{0\}\cup\mathcal{D}^*_{k})=0$ or 1 if $\exists p\nmid d^*$ makes $\nu_{p}(\{0\}\cup\mathcal{D}'^*_{k})=p$, which can't happen to the $k$-tuple jumping champion. On the other hand, from (\ref{c}) and (\ref{d})
\begin{align}
   \mathfrak{S}(\{0\}\cup\lfloor(\log x)^{1/2}\rfloor*\mathcal{K})\frac{x}{\log^{k+1}x}(1-o(1))\le N(x,\mathcal{D}^*_{k})\le\mathfrak{S}(\{0\}\cup\mathcal{D}^*_{k})\frac{x}{\log^{k+1}x}(1+o(1)).\notag
\end{align}
Hence
\begin{align}
\label{h}
   \frac{\mathfrak{S}(\{0\}\cup\mathcal{D}^*_{k})}{\mathfrak{S}(\{0\}\cup\lfloor(\log x)^{1/2}\rfloor*\mathcal{K})}\ge1-o(1).
\end{align}
From (\ref{g}) and (\ref{h}) we obtain
\begin{align}
\label{i}
   \Big(1+\frac{\nu_{p^*}(\{0\}\cup\mathcal{D}'^*_{k})-1} {p^*-\nu_{p^*}(\{0\}\cup\mathcal{D}'^*_{k})}\Big)\le\frac{\mathfrak{S}(\{0\}\cup\lfloor\log^{k+2}x\rfloor*\mathcal{K})} {\mathfrak{S}(\{0\}\cup\lfloor(\log x)^{1/2}\rfloor*\mathcal{K})}(1+o(1)),
\end{align}
while
\begin{align}
   \frac{\mathfrak{S}(\{0\}\cup\lfloor\log^{k+2}x\rfloor*\mathcal{K})} {\mathfrak{S}(\{0\}\cup\lfloor(\log x)^{1/2}\rfloor*\mathcal{K})}\le\sum_{\frac12(1-\epsilon)\log\log x\le p\le(k+2)(1+\epsilon)\log\log x}\frac{p-2}{p-(k+1)}.\notag
\end{align}
Then an argument similar to the deduction of (\ref{f}) from (\ref{e}) gives
\begin{align}
\label{j}
   \frac{\mathfrak{S}(\{0\}\cup\lfloor\log^{k+2}x\rfloor*\mathcal{K})} {\mathfrak{S}(\{0\}\cup\lfloor(\log x)^{1/2}\rfloor*\mathcal{K})}\le1+O\Big(\frac1{\log\log\log x}\Big).
\end{align}
Therefore, from (\ref{i}) and (\ref{j}) we have
\begin{align}
    \Big(1+\frac{\nu_{p^*}(\{0\}\cup\mathcal{D}'^*_{k})-1} {p^*-\nu_{p^*}(\{0\}\cup\mathcal{D}'^*_{k})}\Big)\le(1+o(1))\notag
\end{align}
with $2\le\nu_{p^*}(\{0\}\cup\mathcal{D}'^*_{k})\le\min(k+1, p-1)$. This means that $p^*\rightarrow\infty$ as $x\rightarrow\infty$. Hence we have that any fixed prime $p^*$ must divide every element of a $k$-tuple jumping champions for any given $k\ge2$ and sufficiently large $x$. Then we have proved Theorem \ref{C1}.
\section{proof of theorem \ref{A2}}
From section 4, if $\mathcal{D}^*_k$ is a $k$-tuple jumping champion, it must be that $d^*_k\le\log^{k+1}x$. With the condition of Theorem \ref{A2}, we announce that $d^*_k=o(\log x)$. If not, taking $H=\log x/(\log\log\log x)^{1/2}$ in part (ii) of Lemma \ref{D2}, we have
\begin{align}
    N_k(x,\mathcal{D}^*_{k})\le\mathfrak{S}(\{0\}\cup\mathcal{D}^*_{k})\frac{x}{\log^{k+1}x}\bigg \{1-\frac1{(\log\log\log x)^{1/2}}+o\bigg(\frac1{(\log\log\log x)^{1/2}}\bigg)\bigg\}.\notag
\end{align}
Then using (\ref{f}) with $A=k+2$ and (\ref{j}), we have
\begin{align}
   \mathfrak{S}(\{0\}\cup\mathcal{D}^*_{k})&\le\mathfrak{S}(\{0\}\cup\lfloor(\log x)^{1/2}\rfloor*\mathcal{K})\frac{\mathfrak{S}(\{0\}\cup\mathcal{D}^*_{k})}{\mathfrak{S}(\{0\}\cup \lfloor\log^{k+2}x\rfloor*\mathcal{K})} \frac{\mathfrak{S}(\{0\}\cup\lfloor\log^{k+2}x\rfloor*\mathcal{K})} {\mathfrak{S}(\{0\}\cup\lfloor(\log x)^{1/2}\rfloor*\mathcal{K})}\notag\\
   &\le\mathfrak{S}(\{0\}\cup\lfloor(\log x)^{1/2}\rfloor*\mathcal{K})\bigg \{1+O\bigg(\frac1{\log\log\log x}\bigg)\bigg\}.\notag
\end{align}
It's easy to see that $\lfloor(\log x)^{1/2}\rfloor*k\le k\log^{1/2} x$. Then, from part (i) of Lemma \ref{D2}, we have
\begin{align}
   N_k(x,\mathcal{D}^*_{k})&\le\mathfrak{S}(\{0\}\cup\lfloor(\log x)^{1/2}\rfloor*\mathcal{K})\bigg(1-\frac1{(\log\log\log x)^{1/2}}(1+o(1))\bigg)\notag\\
   &< N_k(x,\lfloor(\log x)^{1/2}\rfloor*\mathcal{K}),\notag
\end{align}
which can't happen to $\mathcal{D}_k^*$. Hence if $\mathcal{D}_k^*$ is a $k$-tuple jumping champion, it must satisfy $d_k^*<H=o(\log x)$.

We also claim that $d^*_k\gg(1-\delta)\frac{\log x}{(\log\log x)^2}$ for any given $\delta>0$. If not, from the famous prime number theorem, we can find prime $p'\le\log\log x$ with $p'\nmid d^*$. It's obvious that $p'd^*_k\le(1-\delta)\Big(\frac{\log x}{\log\log x}\Big)$. Since $\nu_{\{0\}\cup\mathcal{D}^*_{k}}(p')\ge2$, it's easy to see
\begin{align}
   \frac{\mathfrak{S}(\{0\}\cup\mathcal{D}^*_{k})}{\mathfrak{S}(\{0\}\cup p'*\mathcal{D}^*_{k})}&=\bigg(1-\frac{\nu_{\{0\}\cup\mathcal{D}^*_{k}}(p')}{p'}\bigg)\bigg(1-\frac{1}{p'}\bigg)^{-1} \notag\\
   &\le1-\frac1{\log\log x}.\notag
\end{align}
Then, from part (i) of Lemma \ref{D2}, we have
\begin{align}
   N_k(x,\mathcal{D}^*_{k})&\le\mathfrak{S}(\{0\}\cup\mathcal{D}^*_{k})\frac{x}{\log^{k+1}x}\bigg (1+o\Big(\frac1{(\log\log x)^{2}}\Big)\bigg)\notag\\
   &\le\mathfrak{S}(\{0\}\cup p'*\mathcal{D}^*_{k})\frac{x}{\log^{k+1}x}\bigg (1+o\Big(\frac1{(\log\log x)^2}\Big)\bigg)\bigg(1-\frac1{\log\log x}\bigg)\notag\\
   &\le N_k(x,p'*\mathcal{D}^*_{k})\bigg(1-\frac{1-\epsilon}{\log\log x}\bigg)\bigg(1-(1+\epsilon)\frac{p'*d_k^*}{\log x}\bigg)^{-1}\notag\\
   &\le N_k(x,p'*\mathcal{D}^*_{k})(1-\frac\delta2\frac1{\log\log x})\notag\\
   &< N_k(x,p'*\mathcal{D}^*_{k}),\notag
\end{align}
while this can't happen to a $k$-tuple jumping champion. Hence it holds that $d^*_k\ge(1-\delta)\frac{\log x}{(\log\log x)^2}$ for any given $\delta>0$.

We now come to prove $d^*$ is square-free for $(1-\delta)\frac{\log x}{(\log\log x)^2}\le d_k^*=o(\log x)$ with $0<\delta<1$ given. Let $p''$ be a prime that $p''^2\mid d^*$ and $\mathcal{D}^0_k=\frac1{p''}*\mathcal{D}^*_k$. From part (i) of Lemma \ref{D2}, we have
\begin{align}
   N_k(x,\mathcal{D}^*_{k})&=\mathfrak{S}(\{0\}\cup\mathcal{D}^*_{k})\frac{x}{\log^{k+1}x}\bigg\{1-\frac{d^*_{k}}{\log x}+o\bigg(\frac{d^*_{k}}{\log x}\bigg)+o\bigg(\frac1{(\log\log x)^2}\bigg)\bigg\}\notag\\
   &=\mathfrak{S}(\{0\}\cup\mathcal{D}^0_{k})\frac{x}{\log^{k+1}x}\bigg(1-\frac{d^*_{k}}{\log x}+o\Big(\frac{d^*_{k}}{\log x}\Big)\bigg)\notag\\
   &=N_k(x,\mathcal{D}^0_{k})\bigg(1-\frac{d^*_{k}}{\log x}+o\Big(\frac{d^*_{k}}{\log x}\Big)\bigg)\bigg(1-\frac{d^*_{k}}{p''\log x}+o\Big(\frac{d^*_{k}}{\log x}\Big)\bigg)^{-1}\notag\\
   &\le N_k(x,\mathcal{D}^0_{k})\bigg(1-\frac{d^*_{k}}{3\log x}\bigg)\notag\\
   &< N_k(x,\mathcal{D}^0_{k}).\notag
\end{align}
However, this is against the definition of the $k$-tuple jumping champion. Therefore, we have proved that $d^*$ is square-free and obtain Theorem \ref{A2}.
\section{proof of lemma \ref{C2}}
The orginal asymptotic formula of the average of the singular series was given by Gallagher \cite{Gallagher} who proved that
\begin{align}
   \sum_{{1\le d_1,d_2,\cdots d_k\le D}\atop{\text{distinct}}}\mathfrak{S}(\mathcal{D}_k)\sim D^k\notag
\end{align}
In 2004 Montgomery and Soundararajan \cite{Montgomery} proved that, for a fixed $k\ge2$,
\begin{align}
   \sum_{{1\le d_1,d_2,\cdots d_k\le D}\atop{\texttt{distinct}}}\mathfrak{S}(\mathcal{D}_k)= D^k-\bigg({k\atop2}\bigg)D^{k-1}\log D+\bigg({k\atop2}\bigg)(1-\gamma-\log2\pi)D^{k-1}+O(D^{k-3/2+\epsilon}),\notag
\end{align}
where $\gamma$ is Euler's constant. This work strengthens Gallagher's asymptotic formula.

Compared to these formulas which concerned with the average of the singular series over all the components of $\mathcal{D}_k$, in order to determine the precise point of transition between jumping champions, Odlyzko, Rubinstein and Wolf \cite{Odlyzko} proved a asymptotic formulas for the special type of singular series average
\begin{align}
   \sum_{1\le d_1<d_2<\cdots <d_{k-2}< D}\mathfrak{S}(0,d_1,d_2,\cdots,d_{k-2},D)= \mathfrak{S}(D)\frac{D^{k-2}}{(k-2)!}+R_k(D)\notag
\end{align}
with $R_k(D)\ll_k D^{k-2}/\log\log D$. They also presented numerical evidence that suggests that $R_k(D)\ll_k\mathfrak{S}(D)D^{k-3}\log D$. In \cite{Goldston1}, Goldston and Ledoan announced that they have proved $R_k(D)\ll_k D^{k-3+\epsilon}$ for any $\epsilon>0$, but didn't give the proof.

In order to prove the jumping champion conjecture, Goldston and Ledoan, in \cite{Goldston1}, proved the following special type of singular series average, which is different from asymptotic formulas given above,
\begin{align}
   \sum_{1\le d_1<d_2<\cdots <d_{k-2}< H}\mathfrak{S}(0,d_1,d_2,\cdots,d_{k-2},D)= \mathfrak{S}(D)\frac{H^{k-2}}{(k-2)!}(1+o(1))\notag
\end{align}
with $k\ge3$ and $D^\epsilon\le H\le D$. In this paper, we improved this asymptotic formula and actually proved that
\begin{align}
   \sum_{1\le d_1<d_2<\cdots <d_{k-2}< H}\mathfrak{S}(0,d_1,d_2,\cdots,d_{k-2},D)= \mathfrak{S}(D)\frac{H^{k-2}}{(k-2)!}\Big(1+O_k(\frac{D^\epsilon}{H^{1/2}})\Big)\notag
\end{align}
for any $H\le D$. This formula can be deduced easily from Lemma \ref{C2}.

We now come to the proof of Lemma \ref{C2}.

\proof First observe that if $\mathfrak{S}(\mathcal {D}_k)=0$ then $\mathfrak{S}(\mathcal {D}_k\cup{d_0})=0$ and the Lemma holds trivially. Therefore, we assume $\mathfrak{S}(\mathcal {D}_k)\neq0$. Let
\begin{align}
   \mathscr{S}_{d_0}=\frac{\mathfrak{S}(\mathcal {D}_k\cup{d_0})}{\mathfrak{S}(\mathcal {D}_k)}=\prod_p(1+a(p,v_{\mathcal{D}_k\cup \{d_0\}}(p)))\notag,
\end{align}
where
\begin{align}
   a(p,v_{\mathcal{D}_k\cup \{d_0\}}(p))=\frac{(v_{\mathcal {D}_k}(p)-v_{\mathcal {D}_k\cup\{d_0\}}(p)+1)p-v_{\mathcal {D}_k}(p)}{(p-v_{\mathcal {D}_k}(p))(p-1)}.\notag
\end{align}
We now let
\begin{align}
   \Delta_{d_0}=\prod_{1\le i\le k}|d_i-d_0|\notag
\end{align}
and note that
\begin{align}
   v_{\mathcal {D}_k\cup\{d_0\}}(p)=\left\{
\begin{array}
    {r@{\quad:\quad}l}
    v_{\mathcal {D}_k}(p)+1 & p\nmid\Delta_{d_0}; \\
    v_{\mathcal {D}_k}(p)  & p\mid\Delta_{d_0} .
\end{array}
\right.\notag
\end{align}
It follows that
\begin{align}
\label{b3}
   a(p,v_{\mathcal{D}_k\cup \{d_0\}}(p))\ll_k
\left\{
\begin{array}
    {r@{\quad:\quad}l}
    p^{-2} &  p\nmid\Delta_{d_0}; \\
    p^{-1}  & p\mid\Delta_{d_0},
\end{array}
\right.
\end{align}
since $v_{\mathcal{D}_k}(p)\le k$ for any $p$.
Hence the product for $\mathscr{S}_{d_0}$ converges. Defining $a_{d_0}(q)$ for square-free $q$ by
\begin{align}
   a_{d_0}(1)=1,\notag
\end{align}
and
\begin{align}
   a_{d_0}(q)=\prod_{p\mid q}a(p,v_{\mathcal{D}_k\cup \{d_0\}}(p)),\notag
\end{align}
we get
\begin{align}
   \mathscr{S}_{d_0}=\sum_q\mu^2(q)a_{d_0}(q).\notag
\end{align}
It's obvious that the series is convergent.

Let $C$ be a large enough positive constant depending only on $k$. For large $q$, puting $q=q_1q_2$ with $q_1\mid\Delta_{d_0}$ and $(q_2,\Delta_{d_0})=1$, it's obvious that the number of such $q_1$ is $O(h^\epsilon)$. Then we have,
\begin{align}
   \sum_{q>x}\mu^2(q)|a_{d_0}(q)|&\le\sum_{q_1\mid\Delta_{d_0}}\frac{\mu^2(q_1)C^{\omega(q_1)}}{q_1}\sum_{{q_2>x/q_1}\atop {(q_2,\Delta_{d_0})=1}}\frac{\mu^2(q_2)C^{\omega(q_2)}}{q_2^2}\notag\\ &\ll\sum_{q_1\mid\Delta_{d_0}}\frac{1}{q^{1-\epsilon}_1}\sum_{{q_2>x/q_1}\atop {(q_2,\Delta_{d_0})=1}}\frac{1}{q_2^{2-\epsilon}}\ll\sum_{q_1\mid\Delta_{d_0}}\frac{1}{q^{1-\epsilon}_1} \frac{q_1}{x^{1-\epsilon}}\ll(xh)^\epsilon/x,\notag
\end{align}
with the constant depending only on $k$ and $\epsilon$. It follows that
\begin{align}
\label{a}
   \sum_{{1\le d_0\le H}\atop{d_0\notin\mathcal{D}}}\mathscr{S}_{d_0}=\sum_{q\le x}\mu^2(q)\sum_{{1\le d_0\le H}\atop{d_0\notin\mathcal{D}}}a_{d_0}(q)+O(H(xh)^\epsilon/x),
\end{align}
with the constant depending only on $k$ and $\epsilon$.

The inner sum in (\ref{a}) is
\begin{align}
   \sum_v\Big(\prod_{p\mid q}a(p,v(p))\Big)\Big({\sum}'1+O(1)\Big),\notag
\end{align}
where ${\sum}'1$ stands for sum of the number of integer $d_0$ with $1\le d_0\le H$ which, for each prime $p\mid q$, makes $\mathcal{D}_k\cup\{d_0\}$ occupy exactly $v(p)$ residue classes mod $p$; the outer sum is over all "vectors"$=(\cdots, v(p),\cdots)_{p\mid q}$ with components satisfying $v(p)=v_{\mathcal{D}_k}(p)$ or $v_{\mathcal{D}_k}(p)+1$. Here the error term $O(1)$ comes from the ignoring of the condition $d_0\notin\mathcal{D}_k$. The Chinese remainder theorem gives, for $q\le H$ (we choose $x=H^{1/2}\le H$ at last, so this conditions are satisfied.),
\begin{align}
   {\sum}'1=\bigg(\frac Hq+O(1)\bigg)\prod_{p\mid q}f(p,v(p)),\notag
\end{align}
where $f(p,v(p))$ denotes the residue classes of such $d_0$ that makes $v_{\mathcal{D}_k\cup\{d_0\}}(p)=v(p)$. It follows that
\begin{align}
   f(p,v(p))=\left\{
\begin{array}
    {r@{\quad:\quad}l}
    v_{\mathcal {D}_k}(p)     & v(p)=v_{\mathcal {D}_k}(p); \\
    p-v_{\mathcal {D}_k}(p)   & v(p)=v_{\mathcal {D}_k}(p)+1.
\end{array}
\right.\notag
\end{align}

Thus the inner sum in (\ref{a}) is
\begin{align}
   \bigg(\frac Hq\bigg)A(q)+B(q),\notag
\end{align}
with
\begin{align}
   A(q)&=\sum_v\prod_{p\mid q}a(p,v(p))f(p,v(p)),\notag\\
   B(q)&=\sum_v\prod_{p\mid q}|a(p,v(p))|f(p,v(p))+\sum_v\prod_{p\mid q}|a(p,v(p))|.\notag
\end{align}
We have
\begin{align}
   A(q)&=\prod_{p\mid q}\bigg(\sum_{v(p)}a(p,v(p))f(p,v(p))\bigg),\notag\\
   B(q)&=\prod_{p\mid q}\bigg(\sum_{v(p)}|a(p,v(p))|f(p,v(p))\bigg)+\prod_{p\mid q}\bigg(\sum_{v(p)}|a(p,v(p))|\bigg)\notag
\end{align}
From the definition of $a(p,v(p))$ and $f(p,v(p))$,
\begin{align}
   \sum_{v(p)}a(p,v(p))f(p,v(p))&=\frac{p-v_{\mathcal {D}_k}(p)}{(p-v_{\mathcal {D}_k}(p))(p-1)}v_{\mathcal{D}_k}(p)+\frac{-v_{\mathcal {D}_k}(p)}{(p-v_{\mathcal {D}_k}(p))(p-1)}(p-v_{\mathcal {D}_k}(p))\notag\\
   &=0.\notag
\end{align}
Hence, we have that $A(q)=0$ for $q>1$.

Using the bounds (\ref{b3}) for $a(p,v_{\mathcal{D}_k\cup \{d_0\}}(p))$ and the definition of $f(p,v(p))$, we have
\begin{align}
   B(q)\le C^{\omega(q)}.\notag
\end{align}
Employing this into (\ref{a}), it follows that
\begin{align}
   \sum_{{1\le d_0\le H}\atop{d_0\notin\mathcal{D}}_k}\mathscr{S}_{d_0}&=H+O(\sum_{q\le x}C^{\omega(q)})+O(H(xh)^\epsilon/x)\notag\\
   &=H+O(x^{1+\epsilon})+O(H(hx)^\epsilon/x)\notag\\
   &=H+O(H^{1/2}h^\epsilon),\notag
\end{align}
with choosing $x=H^{1/2}$. Then the Lemma follows.

\end{document}